\documentclass [12pt,a4paper,reqno]{amsart} \textwidth 165mm
\usepackage{amsmath,amsthm}
\usepackage{amsfonts,amssymb}
\allowdisplaybreaks

\newcommand{\be}{\begin{equation}}
\newcommand{\ef}{\end{equation}}
\chardef\bslash=`\\ 

\newtheorem*{thm*}{Theorem}

\theoremstyle{definition}

\newtheorem*{remark*}{Remarks}
\newtheorem*{defn*}{Definition}
\theoremstyle{remark} \theoremstyle{remark*}

\newcommand{\G}{\Gamma}
\newcommand{\wt}{\widetilde}
\newcommand{\wh}{\widehat}
\newcommand{\fc}{\frac}

\newcommand{\iy}{\infty}
 \renewcommand{\sectionmark}[1]{}

\newcommand{\Be}{Beltrami}
\newcommand{\hol} {holomorphic}
\newcommand{\qc} {quasiconformal}

\newcommand{\Te} {Teichm\"{u}ller}

 \usepackage{amsfonts}
\newcommand{\field}[1]{\mathbb{#1}}
\newcommand{\g}{\gamma}

\newcommand{\dl}{\delta}
\newcommand{\D}{\field{D}}

\newcommand{\z}{\zeta}
\newcommand{\ov}{\overline}
\newcommand{\vp}{\varphi}
\newcommand{\hC}{\wh{\field{C}}}
\newcommand{\C}{\field{C}}

\newcommand{\B}{\mathbf{B}}
\newcommand{\T}{\mathbf{T}}
\newcommand{\Belt}{\mathbf{Belt}}

\newcommand{\Teich}{\operatorname{Teich}}

\newcommand{\dist}{\operatorname{dist}}

\newcommand{\x} {\mathbf x}

\renewcommand{\a} {\alpha}

\newcommand{\ld}{\lambda}
\newcommand{\kp}{\kappa}

\begin{document}

\title{A new look at Krzyz's conjecture}
\author{Samuel L. Krushkal}

\begin{abstract}  Recently the author has presented a new approach
to solving extremal problems of geometric function theory. It
involves the Bers isomorphism theorem for Teichm\"{u}ller spaces of
punctured Riemann surfaces.

We show here that this approach, combined with quasiconformal
theory, can be also applied to nonvanishing holomorphic functions
from $H^\infty$. In particular this gives a proof of an old open Krzyz conjecture for such functions and of its generalizations.

The unit ball $H_1^\infty$ of $H^\infty$ is naturally embedded into
the universal Teichm\"{u}ller space, and the functions $f \in
H_1^\infty$ are regarded as the Schwarzian derivatives of univalent
functions in the unit disk.
\end{abstract}

\date{\today\hskip4mm({krzyzLook1.tex})}

\maketitle

\bigskip

{\small {\textbf {2010 Mathematics Subject Classification:} Primary:
30C50, 30C55, 30H05; Secondary: 30F60}

\medskip


\bigskip

\markboth{Samuel L. Krushkal}{A new look at Krzyz's conjecture}
\pagestyle{headings}

\bigskip
\centerline{\bf 1. INTRODUCTORY REMARKS AND RESULTS}

\bigskip\noindent
{\bf 1.1. Nonvanishing holomorphic functions and Krzyz's
conjecture}. Consider the circular rings 
$$
\mathcal A_\rho = \{\rho < |z| < 1\}, \quad \rho \ge 0,
$$ 
and denote by $H(\D, \mathcal A_\rho)$ the collection of holomorphic functions $f$ from the unit disk $\D$ into $\mathcal A_\rho$. 
Regarding the points of $\mathcal A_\rho$ as the constant functions on $\D$, one
obtains an embedding of this ring as a subset in the space $H^\iy$
of bounded holomorphic functions in $\D$ with sup-norm.
By $H_1^\infty$ we denote the unit ball of this space. 

The collections $H(\D, \mathcal A_\rho)$ broaden monotonically as $\rho \searrow 0$
giving in the limit the class
$$
H_0^\iy := H(\D, \mathcal A_0) = \bigcup_\rho H(\D, \mathcal A_\rho)
$$
of all nonvanishing $H^\iy$-functions. This class has been actively investigated in geometric function theory from the 1940s, in view of interesting deep features of nonvanishing   (see, e.g., \cite{Go}, \cite{R}).

In 1968, Krzyz \cite{Kz} conjectured that for all functions from
$H_0^\iy$ the following bound
 \begin{equation}\label{1}
|c_n| \le 2/e
\end{equation}
is valid for any $n > 1$, with equality only for the function
$\kp_0(z^n)$ and its pre and post rotations about the origin , where
 \begin{equation}\label{2}
\kp_0(z) = \exp \Bigl(\frac{z -1}{z + 1}
\Bigr) = \fc{1}{e} + \fc{2}{e} z - \fc{2}{3e} z^3 + ... \ .
\end{equation}

This conjecture has been investigated by many 
mathematicians, however it still remains open. The estimate (1) was
established only for some initial coefficients $c_n$ including all
$n \le 5$ (see \cite{HSZ}, \cite{PS}, \cite{Sa}, \cite{Sz}, \cite{Ta}). For developments related to this problem, see, e.g.,
\cite{Ba}, \cite{Ho}, \cite{HSZ}, \cite{Kr1}, \cite{LS},
\cite{MSUV}, \cite{R}, \cite{Su}, \cite{Sz}.

\bigskip\noindent
{\bf 1.2. Main theorem}. Put  
 \begin{equation}\label{3}
\a_{\mathcal A_\rho}
:= \max \{|f^\prime(0)|: \ f \in H(\D, \mathcal A_\rho)\}   
\end{equation}
and take a universal holomorphic covering map $\kappa_\rho: \ \D \to \mathcal A_\rho$, on which this maximal value of $|f^\prime(0)|$ is  attained, i.e., $|\kappa_\rho^\prime(0)| = \a_{\mathcal A_\rho}$. 

Every function $f \in \mathcal H(\D, \mathcal A_\rho)$ admits the factorization
 \begin{equation}\label{4}
 f(z) = \kappa_\rho \circ \wh f(z),
\end{equation}
where $\wh f$ is a holomorphic map of the disk $\D$ into itself (hence, it also belongs to $H_1^\iy$).  

In geometric function theory, such a relation is regarded as a subordination of functions $f$ to $\kappa_\rho$; it has been investigated mostly for univalent covers $\kappa_\rho$. Let 
$$ 
k_\rho(z) = c_0^0 + c_1^0 z + \dots + c_n^0 z^n + \dots, \quad|z| < 1.  
$$ 

The existence of extremal functions maximizing the coefficient $c_n(f)$ on $\mathcal H(\D, \mathcal A_\rho)$ follows from compactness of these classes in the weak topology determined by the locally uniform convergence on $\D$.\footnote{Weak compactness of
$H_0^\infty$ is obtained after adding to this class the function $f(z) \equiv 0$.}

The main result of this paper is the following theorem, 
which implies the proof of Krzyz's conjecture. 

\bigskip\noindent
{\bf Theorem 1}. {\it For all $f \in \ H_0^\iy$ and any $n > 1$,
 \begin{equation}\label{5}
\max |c_n| = \kp_0^\prime(0) = 2/e,
\end{equation}
with equality only for the function $\kp_0(z^n)$ and its
compositions with pre and post rotations about the origin.}

\bigskip
This theorem is extended to the spaces $H(\D, \mathcal A_\rho)$ with sufficiently small $\rho$. 

\bigskip\noindent
{\bf Theorem 2}. {\it There is a number $r_0, \ 0 < r_0 < 1$, such that for any $\rho < \rho_0$ every extremal function $f_0$ maximizing $|c_n|$ on the corresponding class $H(\D, \mathcal A_\rho)$ is of the form 
  \begin{equation}\label{6}
 f_0(z) = \epsilon_2 \kappa_\rho(\epsilon_1 z^n) 
\end{equation} 
with $|\epsilon_1| = |\epsilon_2| = 1$. }

\bigskip 
In the case $\rho > 0$, we do not have an assertion on  uniqueness of the covering map $\kappa_\rho$ on which 
the maximal value (3) is attained.

\bigskip\noindent
{\bf 1.3}. To prove Theorems 1 and 2, we apply a new approach in geometric function theory recently presented  in \cite{Kr2}. It involves the Bers isomorphism theorem for Teichm\"{u}ller spaces of punctured Riemann surfaces.

The unit ball $H_1^\infty$ of $H^\infty$ is naturally embedded into
the universal Teichm\"{u}ller space $\T$, and the functions $f \in H_1^\infty$ are regarded as the Schwarzian derivatives of univalent functions in the unit disk.

In fact, this approach allows one to consider also some more general homogeneous polynomial coefficient functionals than $c_n(f)$.

\bigskip\bigskip
\centerline{\bf 2. DIGRESSION TO TEICHM\"{U}LLER SPACES}

\bigskip
We briefly recall some needed results from Teichm\"{u}ller space theory on spaces  in order to prove Theorem 1; the details
can be found, for example, in \cite{Be}, \cite{GL}.

\bigskip\noindent
{\bf 2.1}.  The universal Teichm\"{u}ller space $\T = \Teich (\D)$
is the space of quasisymmetric homeomorphisms of the unit circle
$\mathbb S^1$ factorized by M\"{o}bius maps;  all Teichm\"{u}ller
spaces have their isometric copies in $\T$.

The canonical complex Banach structure on $\T$ is defined by
factorization of the ball of the Beltrami coefficients (or complex
dilatations)
$$
\Belt(\D)_1 = \{\mu \in L_\iy(\C): \ \mu|\D^* = 0, \ \|\mu\| < 1\},
$$
letting $\mu_1, \mu_2 \in \Belt(\D)_1$ be equivalent if the
corresponding \qc \ maps $w^{\mu_1}, w^{\mu_2}$ (solutions to the
Beltrami equation $\partial_{\ov{z}} w = \mu \partial_z w$ with $\mu
= \mu_1, \mu_2$) coincide on the unit circle $\mathbb S^1 = \partial
\D^*$ (hence, on $\ov{\D^*}$). Such $\mu$ and the corresponding maps
$w^\mu$ are called $\T$-{\it equivalent}. The equivalence classes
$[w^\mu]_\T$ are in one-to-one correspondence with the Schwarzian
derivatives
$$
S_w(z) = \left(\frac{w^{\prime\prime}(z)}{w^\prime(z)}\right)^\prime
- \frac{1}{2} \left(\frac{w^{\prime\prime}(z)}{w^\prime(z)}\right)^2
\quad (w = w^\mu(z), \ \ z \in \D^*).
$$

Note that for each locally univalent function $w(z)$ on a simply
connected hyperbolic domain $D \subset \hC$, its Schwarzian
derivative belongs to the complex Banach space $\B(D)$ of
hyperbolically bounded holomorphic functions on $D$ with the norm
$$
\|\vp\|_\B = \sup_D \ld_D^{-2}(z) |\vp(z)|,
$$
where $\ld_D(z) |dz|$ is the hyperbolic metric on $D$ of Gaussian
curvature $- 4$; hence $\vp(z) = O(z^{-4})$ as $z \to \iy$ if $\iy
\in D$. In particular, for the unit disk,
$$
\ld_\D(z) = 1/(1 - |z|^2).
$$

The space $\B(D)$ is dual to the Bergman space $A_1(D)$, a subspace
of $L_1(D)$ formed by integrable holomorphic functions (quadratic
differentials $\vp(z) dz^2$ on $D$), since every linear functional
$l(\vp)$ on $A_1(D)$ is represented in the form
 \be\label{7}
l(\vp) = \langle \psi, \vp \rangle_D = \iint\limits_D \
\ld_D^{-2}(z) \ov{\psi(z)} \vp(z) dx dy
\end{equation}
with a uniquely determined $\psi \in \B(D)$.

The Schwarzians $S_{w^\mu}(z)$ with $\mu \in \Belt(\D)_1$ range over
a bounded domain in the space $\B = \B(\D^*)$. This domain models
the space $\T$. It lies in the ball $\{\|\vp\|_\B < 6\}$ and
contains the ball $\{\|\vp\|_\B < 2\}$. In this model, the
Teichm\"{u}ller spaces of all hyperbolic Riemann surfaces are
contained in $\T$ as its complex submanifolds.

The factorizing projection
$$
\phi_\T(\mu) = S_{w^\mu}: \ \Belt(\D)_1 \to \T
$$
is a holomorphic map from $L_\iy(\D)$ to $\B$. This map is a split
submersion, which means that $\phi_\T$ has local holomorphic
sections (see, e.g., [GL]).

Note that both equations $S_w = \vp$ and $\partial_{\ov z} w = \mu
\partial_z w$ (on $\D^*$ and $\D$, respectively) determine their
solutions in $\Sigma_\theta$ uniquely, so the values $w^\mu(z_0)$
for any fixed $z_0 \in \C$ and the Taylor coefficients $b_1, b_2,
\dots$ of $w^\mu \in \Sigma_\theta$ depend holomorphically on $\mu
\in \Belt(\D)_1$ and on $S_{w^\mu} \in \T$.

\bigskip\noindent
{\bf 2.2}. The points of Teichm\"{u}ller space $\T_1 =
\Teich(\D_{*})$ of the punctured disk $\D_{*} = \D \setminus \{0\}$
are the classes $[\mu]_{\T_1}$ of $\T_1$-{\it equivalent} \Be \
coefficients $\mu \in \Belt(\D)_1$ so that the corresponding \qc \
automorphisms $w^\mu$ of the unit disk coincide on both boundary
components (unit circle $\mathbb S^1 = \{|z| =1\}$ and the puncture
$z = 0$) and are homotopic on $\D \setminus \{0\}$. This space can
be endowed with a canonical complex structure of a complex Banach
manifold and embedded into $\T$ using uniformization.

Namely, the disk $\D_{*}$ is conformally equivalent to the factor
$\D/\G$, where $\G$ is a cyclic parabolic Fuchsian group acting
discontinuously on $\D$ and $\D^*$. The functions $\mu \in
L_\iy(\D)$ are lifted to $\D$ as the \Be \ $(-1, 1)$-measurable
forms  $\wt \mu d\ov{z}/dz$ in $\D$ with respect to $\G$, i.e., via
$(\wt \mu \circ \g) \ov{\g^\prime}/\g^\prime = \wt \mu, \ \g \in
\G$, forming the Banach space $L_\iy(\D, \G)$.

We extend these $\wt \mu$ by zero to $\D^*$ and consider the unit
ball $\Belt(\D, \G)_1$ of $L_\iy(\D, \G)$. Then the corresponding
Schwarzians $S_{w^{\wt \mu}|\D^*}$ belong to $\T$. Moreover, $\T_1$
is canonically isomorphic to the subspace $\T(\G) = \T \cap \B(\G)$,
where $\B(\G)$ consists of elements $\vp \in \B$ satisfying $(\vp
\circ \g) (\g^\prime)^2 = \vp$ in $\D^*$ for all $\g \in \G$.

Due to the Bers isomorphism theorem, the space $\T_1$ is
biholomorphically isomorphic to the Bers fiber space
$$
\mathcal F(\T) = \{(\phi_\T(\mu), z) \in \T \times \C: \ \mu \in
\Belt(\D)_1, \ z \in w^\mu(\D)\}
$$
over the universal space $\T$ with holomorphic projection $\pi(\psi,
z) = \psi$ (see \cite{Be}).

This fiber space is a bounded hyperbolic domain in $\B \times \C$
and represents the collection of domains $D_\mu = w^\mu(\D)$ as a
holomorphic family over the space $\T$. For every $z \in \D$,  its
orbit $w^\mu(z)$ in $\T_1$ is a holomorphic curve over $\T$.

The indicated isomorphism between $\T_1$ and $\mathcal F(\T)$ is
induced by the inclusion map \linebreak $j: \ \D_{*} \hookrightarrow
\D$ forgetting the puncture at the origin via
 \be\label{8}
\mu \mapsto (S_{w^{\mu_1}}, w^{\mu_1}(0)) \quad \text{with} \ \
\mu_1 = j_{*} \mu := (\mu \circ j_0) \ov{j_0^\prime}/j_0^\prime,
\end{equation}
where $j_0$ is the lift of $j$ to $\D$.

In the line with our goals, we slightly modified the Bers
construction, applying quasiconformal maps $F^\mu$ of $\D_{*}$
admitting conformal extension to $\D^*$ (and accordingly using  the
Beltrami coefficients $\mu$ supported in the disk) (cf. \cite{Kr2}).
These changes are not essential and do not affect the underlying
features of the Bers isomorphism (giving the same space up to a
biholomorphic isomorphism).

The Bers theorem is valid for Teichm\"{u}ller spaces $\T(X_0
\setminus \{x_0\})$ of all punctured hyperbolic Riemann surfaces
$X_0 \setminus \{x_0\}$ and implies that $\T(X_0 \setminus \{x_0\})$
is biholomorphically isomorphic to the Bers fiber space $ \mathcal
F(\T(X_0))$ over $\T(X_0)$.

Note that $\B(\G_0)$ has the same elements as the space $A_1(\D^*,
\G_0)$ of integrable holomorphic forms of degree $- 4$ with norm
$\|\vp\|_{A_1(\D^*, \G_0)} = \iint_{\D^*/\G_0} |\vp(z)| dx dy$; and
similar to (10), every linear functional $l(\vp)$ on $A_1(\D^*,
\G_0)$ is represented in the form
$$
l(\vp) = \langle \psi, \vp \rangle_{\D/\G_0} :=
\iint\limits_{\D^*/\G_0} (1 - |z|^2)^2 \ \ov{\psi(z)} \vp(z) dx dy
$$
with uniquely determined $\psi \in \B(\G_0)$.

Any Teichm\"{u}ller space is a complete metric space with intrinsic
Teichm\"{u}ller metric defined by quasiconformal maps. By the
Royden-Gardiner theorem, this metric equals the hyperbolic Kobayashi
metric determined by the complex structure (see, e.g., \cite{EKK}, \cite{GL}).

\bigskip
We do not use here the finite dimensional Teichm\"{u}ller spaces
corresponding to finitely generated Fuchsian groups.

\bigskip\bigskip
\centerline{\bf 3. PROOF OF THEOREM 1}

\bigskip
We carry out the proof in several stages as a consequence of lemmas.

\bigskip\noindent
$\mathbf{1^0}$. \ Let $G$ be a ring subdomain of the disk $\D$
bounded by the unit circle $S^1 = \partial \D$ and a Jordan curve
$\gamma_G$ separating the origin and $S^1$, and let $H(\D, G)$
denote the subspace of functions $f \in H_1^\iy$ mapping $\D$ into
$G$.

\bigskip
We first establish some results characterizing the structure of such
sets $H(\D, G)$.

\bigskip\noindent
{\bf Lemma 1}. {\it Any set $H(\D, G)$ contains an open path-wise
connective subdomain $H^0(\D, G)$ which is dense in the weak
topology of locally uniform convergence on $\D$, and the universal
holomorphic covering map $\kp_G: \ \D \to G$ extends to a
holomorphic map from $H_1^\iy$ onto $H(\D, G)$ (that is, in
$H^\infty$-norm).}

\bigskip
{\bf Proof}. We precede the proof of this main lemma by two
auxiliary lemmas giving other analytic and geometric features of
sets $H^\iy(\D, G)$.

\bigskip\noindent
{\bf Lemma 2}. {\it (a) \ Every function $f \in H(\D, G)$ admits factorization
 \begin{equation}\label{9}
f(z) = \kp_G \circ \wh f(z),
\end{equation}
where $\wh f$ is a holomorphic map of the disk $\D$ into itself (hence, from $H_1^\iy$).

(b) \ Moreover, the relation generates an $H^\iy$-holomorphic map $\mathbf k_G: \ \wh f \mapsto f$ from 
$H_1^\iy$ onto $H(\D, G)$. }

\bigskip\noindent
{\bf Proof}. (a) Due to a general topological theorem, any map $f: M
\to N$, where $M, N$ are manifolds, can be lifted to a covering
manifold $\wh N$ of $N$, under an appropriate relation between the
fundamental group $\pi_1(M)$ and a normal subgroup of $\pi_1(N)$
defining the covering $\wh N$ (see, e.g, \cite{Ma}). This
construction produces a map $\wh f: M \to \wh N$ satisfying
 \begin{equation}\label{10}
f = p \circ \wh f,
\end{equation}
where $p$ is a projection $\wh N \to N$. The map $\wh f$ is
determined up to composition with the covering transformations of
$\wh N$ over $N$ or equivalently, up to choosing a preimage of a
fixed point $x_0 \in \wh N$ in its fiber $p^{-1}(x_0)$. For \hol \
maps and manifolds, the lifted map is also holomorphic.

In our special case, $\kp_G$ is a holomorphic universal covering map
$\D \to G$, and the representation (10) provides the equality (9)
with the corresponding $\wh f$ determined up to covering
transformations of the unit disk compatible with the covering map
$\kp_G$.

For a fixed $z \in \D$, each coefficient $c_n$, of $f$ (and hence
$f$ itself) is a holomorphic function (polynomial) of the initial
coefficients $\wh c_0, \wh c_1, \dots, \wh c_n$ of cover $\wh f$.
Holomorphy in the $H^\iy$ norm stated by the assertion {\it (b)} is a
consequence of a well-known property of bounded holomorphic functions in Banach spaces with sup norm given by the following lemma of Earle \cite{Ea}.

\bigskip\noindent
{\bf Lemma 3}. {\it Let $E, \ T$ be open subsets of complex Banach
spaces $X, Y$ and $B(E)$ be a Banach space of \hol \ functions on
$E$ with sup norm. If $\vp(x, t)$ is a bounded map $E \times T \to
B(E)$ such that $t \mapsto \vp(x, t)$ is holomorphic for each $x \in E$,
then the map $\vp$ is holomorphic. }

\bigskip
Holomorphy of $\vp(x, t)$ in $t$ for fixed $x$ implies the existence
of complex directional derivatives
$$
\vp_t^\prime(x,t) = \lim\limits_{\z\to 0} \fc{\vp(x, t + \z v) -
\vp(x, t)}{\z} = \fc{1}{2 \pi i} \int\limits_{|\xi|=1} \fc{\vp(x, t
+ \xi v)}{\xi^2} d \xi,
$$
while the boundedness of $\vp$ in sup norm provides the uniform
estimate
$$
\|\vp(x, t + c \z v) - \vp(x, t) - \vp_t^\prime(x,t) c v\|_{B(E)}
\le M |c|^2,
$$
for sufficiently small $|c|$ and $\|v\|_Y$.

The map $\mathbf k_\rho: \ \wh f \mapsto f$ is bounded on the ball
$H_1^\iy$. Applying Hartog's theorem on separate holomorphy to the
sums $g(z, t) = \wh f(z) + t \wh h(z)$ of $\wh f \in H_1^\iy, \ \wh
h \in H_1$ and $t$ from a region $B \subset \hC$ so that $g(z, t)
\in H_1^\iy$, one obtains that $g(z, t)$ are jointly holomorphic in
both variables $(z, t) \in \D \times B$. Thus the  restriction of
the map $\mathbf k_G$ onto intersection of the ball $H_1^\iy$ with
any complex line $L = \{\wh f + t \wh h\}$ is $H^\iy$-\hol, and
hence this map is \hol \ as the map $H_1^\iy \to H(\D, G)$, which
completes the proof of Lemma 3.

\bigskip
One can also show that the restriction of the extended map $\mathbf k_G$ to any holomorphic disk
$$
\D(\wh f) = \{t \wh f/\|\wh f\|_\iy: \ |t| < 1\}, \quad \wh f \in H_1^\iy,
$$
is a complex geodesic (cf. \cite{Ve}), hence a local hyperbolic
isometry (preserving such property of the original map $k_G$). We
will not use this fact and therefore do not present here its proof.

\bigskip
Now consider the domains $G \Subset \D$. Fix a point $w_0
\in G$ and take for a decreasing sequence $r_n \to 0$ the connected components $G_n^0$ of widening open sets
$$
G_n = \{w \in G: \ \dist(w, \partial G) > r_n\}, \ \ G_n \Subset
G_{n+1}, \quad n = 1, 2, \dots \ , 
$$
exhausting this domain and containing $w_0$. Let
$$
H^0(\D, G) = \bigcup_n H(\D, G_n).
$$
The set $H^0(\D, G)$ is open and contains, in particular, all functions $f \in H(\D, G)$ holomorphic on the closed disk $\ov{\D}$.

\bigskip\noindent
{\bf Lemma 4}. {\it (a) \ For any fixed $n$, every function $f \in H(\D, G_n)$ continuous in the the closed disk $\ov{\D}$ has a neighborhood $U(f, \epsilon_n)$ in $H_1^\iy$, which contains only
the functions belonging to $H^0(\D, G)$.

(b) \ Each of the sets $ H(\D, G_n)$ and $H^0(\D, G)$ is path-wise connective in $H_1^\iy$; therefore, the union $H^0(\D, G)$ is a domain in $H_1^\iy$.   }

\bigskip \noindent
{\bf Proof}. To prove the assertion {\it (a)}, assume to  the contrary that, for some $n = n_0$, such a number $m(n_0)$ does not exist. Then there is a function $f_0 \in H(\D, G_{n_0})$, a sequence of functions $f_m \in H(\D, G_{n_0})$ convergent to $f_0$ so that
 \be\label{11}
\lim\limits_{m\to \iy} \|f_m - f_0\|_{H^\iy} = 0,
\end{equation}
and a sequence of points $z_m \in \D$  convergent to $z_0 \in \D$, for which we will have $f_m (z_m) \in G_m$ for all $m \ge n_0$, and
 \be\label{12}
\lim\limits_{m\to \infty} f_m(z_m) = a \in \D \setminus \ov{G}.
\end{equation}

We approximate $f_m(z)$ by functions $f_{m,r}(z) = f_m(r_m z)$
(holomorphic in $\ov{\D}$), taking $r_m$ so close to $1$ that the
equality (11) is preserved for $f_{m,r}$. Then the uniform
convergence of $f_m$ and $f_{m,r}$ to $f_0$ on compact subsets of
$\D$ immediately implies that the limit $a$ in (15) must be equal to $f_0(z_0)$, and therefore it belongs to $G_{n_0}$. This proves part {\it (a)}.

\bigskip
To show that each $H(\D, G_n)$ is path-wise connective,  take its arbitrary distinct points $f_1, \ f_2$. By (9),
$$
f_j = \kp_{G_n} \circ \wt f_j \ (\wt f_j \in H_1^\iy), \quad j = 1,
2.
$$
Connecting the covers $\wt f_1$ and $ \wt f_2$ in $H_1^\iy$ by the
line interval $l_{1,2}(t) = t \wt f_1 + (1 - t)\wt f_2 \ 0 \le t \le
1$, one obtains a path $\kp_{G_n} \circ l_{1,2}: \ [0, 1] \to H(\D,
G_n)$ connecting $f_1$ with $f_2$, completing the proof of Lemmas 4 and 1.

\bigskip
Observe that Lemma 4 does not contradict to existence for $f_0 \in H(\D, G)$ of sequences $\{f_n\} \in H_1^\iy$ convergent to $f_0$
only locally uniformly in $\D$ and taking some values in $\D \setminus G$.

\bigskip
Using the homotopy $\wh f_t(z) = t \wh f(tz)$ of the cover functions
and representation (9), one concludes that the domain $H^0(\D, G)$
is dense in the set $H(\D, G)$ in the weak topology. Hence,
$$
\sup_{H^0(\D,G)} |J(f)| = \max_{H(\D,G)} |J(f)|
$$
for any holomorphic functional $J(f)$. This follows also from the fact that all $f \in H(\D, G)$ holomorphic on the closure $\ov G$ of domain $G$ belong to $H^0(\D, G)$.

\bigskip
Note also that for $G = \mathcal A_\rho$, the distinguished domain
$H^0(\D, \mathcal A_\rho)$ preserves circular symmetry, i.e., it contains the nonvanishing functions $f \in H_1^*$ together with their compositions with pre and post rotations about the origin.

\bigskip\noindent
$\mathbf 2^0$. \ In the case of the punctured disk $\D_{*}= \mathcal A_0$ Lemma 4 admits some strengthening.

\bigskip\noindent
{\bf Lemma 5}. {\it Each point $f \in H_0^\iy$ has a
neighborhood (ball) $U(f, \epsilon)$ in $H^\iy$, which entirely
belongs to $\mathcal B$, i.e., contains only nonvanishing functions
on the disk $\D$. Take the maximal balls $U(f, \epsilon)$ with such
property. Then their union
$$
U_0 = \bigcup_{f\in H^\iy} U(f, \epsilon)
$$
is a domain in the space $H^\iy$. }

\bigskip\noindent
{\bf Proof}. {\it Openness}: It suffices to show that for each $r > 1$ and $r^\prime < r $, every function $f \in H_0^\iy$ has a neighborhood $U(f, \epsilon(r))$ in $H^\iy(\D_{r^\prime})$, which contains only nonvanishing functions on $\D_{r^\prime} = \{|z| < r^\prime\}$. For $r^\prime = 1$, this gives the first assertion of
the lemma.

Assume the contrary. Then (for some $r > 1$ and $r^\prime < r$)
there exist a function $f_0 \in \mathcal B_r$ and the sequences of
functions $f_n \in H^\iy(\D_{r^\prime})$ convergent to $f_0$,
 \be\label{13}
\lim\limits_{n\to \iy} \|f_n - f_0\|_{H^\iy(\D_{r^\prime})} = 0
\end{equation}
and of points $z_n \in \D$ convergent to $z_0, \ |z_0| \le r^\prime$
such that $f_n (z_n) = 0 \ (n = 1, 2, \dots)$.

In the case $|z_0| < r^\prime$, we immediately reach a contradiction,
because then the uniform convergence of $f_n$ on compact sets in
$\D_{r^\prime}$ implies $f_0(z_0) = 0$, which is impossible.

The case $|z_0| = r^\prime$ requires other arguments. Since $f_0$ is
\hol \ and does not vanish on the closed disk $\ov{\D_r}$,
$$
\min_{|z| \le r^\prime} |f_0(z)| = a > 0.
$$
Hence, for each $z_n$,
$$
|f_n(z_n) - f_0(z_n)| = |f_0(z_n)| \ge a,
$$
and by continuity, there exists a neighborhood $\D(z_n, \dl_n) =
\{|z - z_n < \dl_n\}$ of $z_n$ in $\D_{r^\prime}$, in which $|f_n(z)
- f_0(z)| \ge a/2$ for all $z$. This implies
$$
 \|f_n - f_0\|_{H^\iy(\D_{r^\prime})} \ge
\max_{\D(z_n, \dl_n)} |f_n(z) - f_0(z)|
 \ge \fc{a}{2}.
 $$
This inequality must hold for all $n$, contradicting (13).

The {\it connectedness} of the union $\mathcal H$ is established similar to Lemma 4.

\bigskip\noindent
$\mathbf 3^0$. \ The next step in the proof of Theorem 1 is to construct a holomorphic embedding of the unit ball $H_1^\infty$ into some Teichm\"{u}ller spaces.

Any function $g$ from this ball $H_1^\iy$ belongs to the space $\B =
\B(\D)$ of hyperbolically bounded holomorphic functions $f(z)$
regarding as holomorphic quadratic differentials $f(z) dz^2$ on the
unit disk, with norm
$$
\|f\|_\B = \sup_\D(1 - |z|^2)^2 |f(z)| < 1. 
$$
Hence, such $f$ is the Schwarzian derivative $S_w$ of a univalent function $w(z)$ in the unit disk $\D$ solving the differential equation $S_w = f$. This $w$ is determined up to a Moebius map of the sphere $\hC$.

This implies a holomorphic embedding $\iota$ of the ball
$H_1^\infty$ into the universal Teichm\"{u}ller space.

To determine $w$ uniquely (and ensure the holomorphic dependence of
$w$ from $S_w$), we shall use the following normalization.

Consider similar to \cite{Kr2} the family $\wh S(1)$ of univalent
functions on $\D$ which is the completion in the topology of locally
uniform convergence on $\D$ of the set of univalent functions $w(z)
= a_1 z + a_2 z^2 + \dots$ with $|a_1| = 1$, having quasiconformal
extensions across the unit circle $\mathbb S^1 =
\partial \D$ to the whole sphere $\hC = \C \cup \{\iy\}$ which
satisfy $w(1) = 1$.

Equivalently, this family is a disjunct union
$$
\wh S(1) = \bigcup_{- \pi \le \theta < \pi} S_\theta(1),
$$
where $S_\theta(1)$ consists of univalent functions $w(z) = e^{i
\theta} z + a_2 z^2 + \dots$ with quasiconformal extensions to $\hC$
satisfying $w(1) = 1$ (also completed in the indicated weak
topology).

This family is closely related to the canonical class $S$ of
univalent functions $w(z)$ on $\D$ normalized by $w(0) = 0, \ w^\prime(0) = 1$. Every $w \in S$ has its representative $\wh w$ in
$\wh S(1)$ (not necessarily unique) obtained by pre and post compositions of $w$ with rotations $z \mapsto e^{i \alpha} z$ about
the origin, related by 
  \be\label{14}
w_{\tau, \theta}(z) =  e^{- i \theta} w(e^{i \tau} z) \quad
\text{with} \ \ \tau = \arg z_0,
\end{equation}
where $z_0$ is a point for which $w(z_0) = e^{i \theta}$ is a common point of the unit circle and the the boundary of domain $(\D)$. In the general case, the equality $w(1) = 1$  in terms of the Carath\'{e}odory prime ends.

The relation (14) implies, in particular, that the functions
conformal in the closed disk $\ov \D$ are dense in each class
$S_\theta(1)$. Such a dense subset is formed, for example, by the
images of the homotopy functions $[f]_r(z) = \fc{1}{r} f(r z)$ with
real $r \in (0, 1)$.

The inverted functions $F(z) = 1/f(1/z)$ for $f \in \wh S(1)$ form
the corresponding classes $\Sigma_\theta(1)$ of nonvanishing
univalent functions on the complementary disk
$$
\D^* = \{z \in \hC: \ \ |z| > 1\}
$$
with expansions
$$
F(z) =  e^{- i \theta} z + b_0 + b_1 z^{-1} + b_2 z^{-2} + \dots,
\quad  F(1) = 1,
$$
and $\wh \Sigma(1) = \bigcup_\theta \Sigma_\theta(1)$.

Simple computations yield that the coefficients $a_n$ of $f \in
S_\theta(1)$ and the corresponding coefficients $b_j$ of $F(z) =
1/f(1/z) \in \Sigma_\theta(1)$ are related by
$$
b_0 + e^{2i \theta} a_2 = 0, \quad b_n + \sum \limits_{j=1}^{n}
\epsilon_{n,j}  b_{n-j} a_{j+1} + \epsilon_{n+2,0} a_{n+2} = 0,
\quad n = 1, 2, ... \ ,
$$
where $\epsilon_{n,j}$ are the entire powers of $e^{i \theta}$. This
successively implies the representations of $a_n$ by $b_j$ via
 \be\label{15}
a_n = (- 1)^{n-1} \epsilon_{n-1,0}  b_0^{n-1} - (- 1)^{n-1} (n - 2)
\epsilon_{1,n-3} b_1 b_0^{n-3} + \text{lower terms with respect to}
\ b_0.
\end{equation}

The coefficients $\a_n$ of Schwarzians
$$
S_w(z) = \sum_0^\infty \a_n z^n
$$
are represented as polynomials of $n + 2$ initial coefficients of $w \in S_\theta(1)$ and, in view of (15), as polynomials of $n + 1$ initial coefficients of the corresponding $W \in \Sigma_\theta(1)$;
denote these polynomials by $J_n(w)$ and $\wt J_n(W)$, respectively.
These polynomial functionals are naturally extended to the whole
classes $\wh S(1)$ and $\wh \Sigma(1)$.

\bigskip\noindent
$\mathbf{4^0}$. \ As was mentioned above, any $f \in H^\iy$ belongs
to the space $\B$ and
$$
\|f\|_\B = \sup_\D(1 - |z|^2)^2 |f(z)| < \|f\|_\iy.
$$
By the well-known Ahlfors-Weill theorem, any $g \in \B$ with norm
$\|g\|_\B  = k < 2$ is the Schwarzian derivative $S_w = g$ of a
function $w$ which is univalent on the disk $\D$ and admits
$k$-quasiconformal extension across the unit circle $\{|z| =1\}$ to
$\hC$ with Beltrami coefficients
$$
\nu_{S_w}(\z) = \partial_{\ov{\z}} w/\partial_\z w = - \fc{1}{2}
(|\z|^2 - 1)^2 \fc{\z^2}{\ov{\z}^2} S_w\Bigl(\fc{1}{\ov \z}\Bigr).
$$
We shall denote this holomorphic embedding of the ball $H_1^\infty$
into the space $\T$ modeled by Schwarzians in $\D^*$ again by
$\iota$. The image of $H^\iy$ under this embedding is a noncomplete
linear subspace in $\B$ so that $\iota H_1^\iy$ is a complex subset
of the unit ball in $\B$, and the image of the distinguished domain
$H^0(\D, \mathcal A_\rho)$ is a complex submanifold in $\T$.

Another important property of the set $\iota H^0(\D, \mathcal
A_\rho)$ is given by the following two lemmas.

\bigskip\noindent
{\bf Lemma 6}. {\it Let $f(z) = \sum\limits_0^\iy c_n z^n \in H^\iy$
and $s_m(z) = \sum\limits_0^{m-1} c_n z^n$. Then }
 \be\label{16}
\lim\limits_{m \to \iy} \|s_m - f\|_\B = 0.
\end{equation}

\bigskip\noindent
{\bf Proof}. It suffices to consider the functions $f$ from the ball
$\{\|f\|_\B < 1/2\}$. Their coefficients $c_j$ are estimated by
$|c_n| < 1/2$ for all $n \ge 0$. Hence,
$$
|s_m(z)| < \fc{1}{2}\sum\limits_0^{m-1} |z|^n < \fc{1}{2(1 - |z|)},
$$
and
$$
\|s_m\|_\B < \frac{1}{2} \ \sup_\D (1 - |z|)(1 + |z|)^2 < 2,
$$
which means that any partial sum $s_n$ for such $f$ lies in the ball
$\{g \in \B: \ \|g\|_\B < 2\}$ and therefore it also belongs to the
space $\T$. Further,
$$
\begin{aligned}
|s_m(z) - f(z)| &= |c_m z^m + c_{m+1} z^{m+1} + \dots| < \fc{1}{2} \
(|z|^m + |z|^{m+1} + \dots)  \\
&= \fc{1}{2} \ \frac{|z|^m}{1 - |z|} < \fc{1}{2^{m+1}} \ \fc{1}{1 -
|z|},
\end{aligned}
$$
which implies
$$
\|s_m(z) - f(z)\|_\B < \fc{1}{2^{m+1}} \ \sup_\D (1 - |z|)(1 +
|z|)^2 < \fc{1}{2^{m`-1}},
$$
and (19) follows, proving the lemma.

\bigskip
Lemmas 1 and 5 imply

\bigskip\noindent
{\bf Lemma 7}. {\it Any point $f(z) = \sum_0^\iy c_n z^n$ from the set $\iota H^0(\D, \mathcal A_\rho)$ in $\T$ is approximated in the 
$\B$-norm by polynomials $s_m (z) = c_0 + c_1 z + \dots + c_m z^m$ with $m \ge m_0(f)$, which also belong to $\iota H^0(\D, \mathcal A_\rho)$ and hence do not have zeros in the disk $\D$. }

\bigskip\noindent
$\mathbf{5^0}$. \ Our next step is to lift both  polynomial functionals $J_n(w)$ and $\wh J_n(W)$ (equivalently $c_n$) onto the Teichm\"{u}ller space $\T_1$. Letting
 \be\label{17}
\wh J_n(\mu) = \wt J_n(W^\mu),
\end{equation}
we lift these functionals  from the sets $S_\theta(1)$ and
$\Sigma_\theta(1)$ onto the ball $\Belt(\D)_1$. Then, under the
indicated $\T_1$-equivalence, i.e., by the quotient map
$$
\phi_{\T_1}: \ \Belt(\D)_1 \to \T_1, \quad \mu \to [\mu]_{\T_1},
$$
the functional $\wt J_n(W^\mu)$ is pushed down to a bounded
holomorphic functional $\mathcal J_n$ on the space $\T_1$ with the same range domain.

Equivalently, one can apply the quotient map $\Belt(\D)_1 \to \T$
(i.e., $\T$-equivalence) and compose  the descended functional on
$\T$ with the natural holomorphic map $\iota_1: \ \T_1 \to \T$
generated by the inclusion $\D_{*} \hookrightarrow \D$ forgetting
the puncture. Note that since the coefficients $b_0, \ b_1, \dots$
of $W^\mu \in \Sigma_\theta$   are uniquely determined by its
Schwarzian $S_{W^\mu}$, the values of $\mathcal J_n$ in the points
$X_1, \ X_2 \in \T_1$ with $\iota_1(X_1) = \iota_1(X_2)$ are equal.

Now, using the Bers isomorphism theorem, we regard the points of the
space $\T_1$ as the pairs $X_{W^\mu} = (S_{W^\mu}, W^\mu(0))$, where
$\mu \in \Belt(\D)_1$ obey $\T_1$-equivalence (hence, also
$\T$-equivalence). Denote (for simplicity of notations) the
composition of $\mathcal J_n$ with biholomorphism $\T_1 \cong \mathcal
F(\T)$ again by $\mathcal J_n$. In view of (8) and (17), it is
presented on the fiber space $\mathcal F(\T)$ by
 \be\label{18}
\mathcal J_n(X_{W^\mu}) = \mathcal J_n(S_{W^\mu}, \ t), \quad t = W^\mu(0).
\end{equation}
This yields a logarithmically plurisubharmonic functional $|\mathcal
J_n(S_{W^\mu}, t)|$ on $\mathcal F(\T)$.

Note that since the coefficients $b_0, \ b_1, \dots$ of $W^\mu \in
\Sigma_\theta$   are uniquely determined by its Schwarzian
$S_{W^\mu}$, the values of $\mathcal J_n$ in the points $X_1, \ X_2
\in \T_1$ with $\iota_1(X_1) = \iota_1(X_2)$ are equal.

\bigskip
We have to estimate a smaller plurisubharmonic functional arising
after restriction of $\mathcal J_n(S_{F^\mu}, \ t)$ onto the the
images in these spaces of the distinguished convex set $\iota
H^0(\D, \mathcal A_\rho)$,  i.e., the functional (17) on the set of
$S_{W^\mu} \in \iota H^0(\D, \mathcal A_\rho)$ and corresponding
values of $t = W^\mu(0)$ which runs over some subdomain
$D_{\rho,\theta}$ in the disk $\{|t| < 4\}$.

We denote this restricted functional by $\mathcal J_{n,0}(S_{W^\mu}, \
W^\mu(0))$ and define in domain $D_{\rho,\theta}$ the function
 \be\label{19}
u_\theta(t) = \sup_{S_{W^\mu}} |\mathcal J_{n,0}(S_{W^\mu}, t)|,
\end{equation}
where the supremum is taken over all $S_{F^\mu} \in \iota H^0(\D,
\mathcal A_\rho)$ admissible for a given $t = W^\mu(0) \in
D_{\rho,\theta}$, that means over the pairs $(S_{W^\mu}, t) \in
\mathcal F(\T)$ with $S_{F^\mu} \in \iota H^0(\D, \mathcal A_\rho)$
and a fixed $t$.

Our goal is to establish that this function inherits the
subharmonicity of $\mathcal J$. This is given by the following basic
lemma.

\bigskip\noindent
{\bf Lemma 8}. {\it The function $u_\theta(t)$ is subharmonic in the
domain $D_{\rho,\theta}$. }

\bigskip\noindent
{\bf Proof}. Consider in the set $\iota H^0(\D, \mathcal A_\rho)$
its $m$-dimensional analytic subsets $V_m$ corresponding to the
partial sums $s_m$ of functions $f \in H^0(\D, \mathcal A_\rho)$
(with $m \ge m_0(f)$). Given such $f$, we define
$$
F(z) = f(1/z)/z^4
$$
and take a univalent solution $W \in \Sigma_\theta$ of the
Schwarzian equation $S_W(z) = F(z)$ on $D^*$. Let $W^\mu$ be one of
its quasiconformal extensions onto $\D$.

Let $W_m$ and $W_m^{\mu_m}$ be the corresponding functions defined
similarly by the partial sums $s_m$ of $f, \ m \ge m_0(f)$. Then the
domains $W_m(\D^*)$ and $W_m^{\mu_m}(\D)$ approximate $W(\D^*)$ and
$W^\mu(\D)$ uniformly (in the spherical metric on $\hC$), and the
points $W_n^{\mu_m}(0)$ are close to $W^\mu(0)$.

One can replace the extensions $W_n^{\mu_m}$ by $\omega_m \circ
W_n^{\mu_m}$, where $\omega_m$ is the extremal quasiconformal
automorphism of domain $W_m^{\mu_m}(\D)$ moving the point
$W_m^{\mu_m}(0)$ into $W^\mu(0)$ and identical on the boundary of
$W_m^{\mu_m}(\D)$ (cf. \cite{Te}). This provides for a prescribed $t
= W^\mu(0)$ the points $S_{W_m^{\mu_m}} \in \mathcal F(\T)$
corresponding to given $s_m \in V_m$.

Now, maximizing the function $\log |\mathcal J_{n,0}(S_{W_m^{\mu_m}},
t)|$ over the manifold $V_m$, i.e., over $S_{W_m^{\mu_m}}$ (with appropriate $m$), one obtains a logarithmically plurisubharmonic function
$$
u_m(t) = \sup_{V_m} |\mathcal J_{n,0}(S_{W_m^{\mu_m}}, t)|, \quad t =
W^\mu(0),
$$
in the domain $D_{\rho,\theta}$ indicated above. We take its upper semicontinuous regularization
$$
u_m(t) = \limsup\limits_{t^\prime \to t} u_m(t^\prime)
$$
(denoted, by abuse of notation, by the same letter as the original
function). The general properties of subharmonic functions in the
Euclidean spaces imply that such a regularization also is
logarithmically subharmonic in $D_{\rho,\theta}$.

In a similar way, taking the limit
$$
u(t) = \limsup\limits_{m \to \infty} u_m(t)
$$
followed by its upper semicontinuous regularization, one obtains a
logarithmically subharmonic function on the domain $D_{\rho,\theta}$. Lemmas 6 and 7 imply that this function coincides with function (19).

\bigskip\noindent
$\mathbf{6^0}$. \ Assume now that $\rho = 0$, hence $\mathcal A_\rho = \D_{*} = \D \setminus \{0\}$. 

We have to establish the value domain of $W^\mu(0)$ for $W^\mu$ running over $\iota H^0(\D, \mathcal A_0)$. 

First, we apply the following generalization of the above construction. Taking a dense countable subset
$$
\Theta =  \{ \theta_1, \theta_2, \dots, \theta _m, \dots \} \subset [-\pi, \pi],
$$
consider the increasing unions of the quotient spaces
 \be\label{20}
\mathcal T_m = \bigcup_{j=1}^m \ \wh \Sigma_{\theta_j}^0/\thicksim \
= \bigcup_{j=1}^m \{(S_{W_{\theta_j}}, W_\theta^\mu(0)) \} \ \simeq
\T_1 \cup \dots \cup \T_1,
\end{equation}
where the equivalence relation $\thicksim$ means $\T_1$-equivalence
on a dense subset $\wh \Sigma^0(1)$ in the union $\wh \Sigma(1)$
formed by all univalent functions $W_\theta(z) = e^{-i \theta_j} z +
b_0 + b_1 z^{-2} + \dots$ on $\D^*$ (preserving $z = 1$) with quasiconformal extension to $\hC$, and
$$
\mathbf W_\theta^\mu(0) := (W_{\theta_1}^{\mu_1}(0), \dots ,
W_{\theta_m}^{\mu_m}(0)).
$$
The Beltrami coefficients  $\mu_j \in \Belt(\D)_1$ are chosen here independently. The corresponding collection $\beta = (\beta_1, \dots, \beta_m)$ of the Bers isomorphisms
$$
\beta_j: \ \{(S_{W_{\theta_j}}, W_{\theta_j}^{\mu_j}(0))\} \to
\mathcal F(\T)
$$
determines a holomorphic surjection of the space $\mathcal T_m$ onto $\mathcal F(\T)$.

Taking also in each union (20) the corresponding collection $\iota_m H_0^\iy$ covering $H^0(\D, \mathcal A_0)$, one obtains in a similar manner to the above the maximal function
 \be\label{21}
u(t) = \sup_\Theta u_{\theta_m}(t) = \sup \{|\mathcal J_{n,0}(S_{W_\theta^\mu}, t)|: \ \theta \in \bigcup_m \iota_m H_0^\iy\}.
\end{equation}
It is defined and subharmonic in domain
$$
D_\rho = \bigcup_{\Theta} D_{\rho,\theta_m}.
$$

Noting that the union of spaces $\mathcal T_m$ possesses  the circular
symmetry inherited from the class $\wh \Sigma(1)$, which is preserved under rotations (14),  
one concludes that this broad domain $D_0$ must be a disk $\D_{r_0} = \{|t| < r_0\}$. 

Now we show that in the case of nonvanishing $H^\iy$ functions this radius $r_0$ is naturally connected with the function (2). This requires a covering estimate of Koebe's type.  

Let $G$ be a domain in a complex Banach space $X = \{\mathbf x\}$
and $\chi$ be a holomorphic map from $G$ into the universal
Teichm\"{u}ller space $\T$ modeled as a bounded subdomain of $\B$.
Consider in the unit disk the corresponding Schwarzian differential
equations
 \be\label{22}
S_w(z) = \chi(\x)
\end{equation}
and pick their univalent solutions $w(z)$ satisfying $w(0) =
w^\prime(0) - 1 = 0$ (hence $w(z) = z  + \sum_2^\infty a_n z^n$).
Put
 \be\label{23}
|a_2^0| = \sup \{ |a_2|: \ S_w \in \chi(G)\},
\end{equation}
and let $w_0(z) = z + a_2^0 z^2 + \dots$ be one of the maximizing
functions.

\bigskip\noindent
{\bf Lemma 10}. {\it (a) For every indicated solution $w(z) = z + a_2 + \dots$ of (22), the image domain $w(\D)$ covers entirely the disk
$\{|w| < 1/(2 |a_2^0|)\}$.

The radius value $1/(2 |a_2^0|)$ is sharp for this collection of
functions, and the circle $\{|w| = 1/(2 |a_2^0|)$ contains points
not belonging to $w(\D)$ if and only if $|a_2| = |a_2^0|$ (i.e., when $w$ is one of the maximizing functions).

(b) The inverted functions
$$
W(\zeta) = 1/w(1/\zeta) = \zeta - a_2^0 + b_1 \zeta^{-1} + b_2 z^{-2} + \dots
$$
map the disk $\D^*$ onto a domain whose boundary is entirely
contained in the disk} $\{|W + a_2^0| \le |a_2^0|\}$.

\bigskip\noindent
{\bf The proof} follows the classical lines of Koebe's $1/4$ theorem (cf. \cite{Go}).

{\it (a)} Suppose that the point $w = c$ does not belong to the image of $\D$ under the map $w(z)$ defined above. Then $c \ne 0$, and the function
$$
w_1(z) = c w(z)/(c - w(z)) = z + (a_2 + 1/c) z^2 + \dots
$$
also belongs to this class, and hence by (23), $|a_2 +1/c| \le |a_2^0|$, which implies
$$
|c| \ge 1/(2 |a_2^0|).
$$
The equality holds only when
$$
|a_2 + 1/c| = |1/c| - |a_2| = |a_2^0| \quad \text{and} \ \ |a_2| = |a_2^0|.
$$

{\it(b)} If a point $\zeta = c$ does not belong to the image $W(\D^*)$, then the function 
$$
W_1(z) = 1/[W(1/z) - c] = z + (c + a_2) z^2 + \dots 
$$
is holomorphic and univalent in the disk $\D$, and therefore,  $|c + a_2| \le |a_2^0|$. The lemma follows.

\bigskip
This lemma implies that the boundary of the range domain of $W^\mu(0)$ is contained in the disk
 \be\label{24}
\D_{2|a_2^0|} = \{W: \ |W| \le 2 |a_2^0|\},
\end{equation} 
and, consequently, $r_0 = 2 |a_2^0|$ and touches from inside the circle $\{|W| = 2 |a_2^0|\}$ at the points corresponding to extremal functions $W_0$ maximizing $|a_2|$ on the closure of the domain $\iota H_0^*$. 

Generically, the extremal value $2 |a_2^0|$ of the radius of covered disk can be attained on several functions $W_0$. 

\bigskip\noindent
$\mathbf{7^0}$. \ We now establish that 
  \be\label{25}
S_{W_0}(z) = \kappa_0(z). 
\end{equation} 
In view of Lemma 1, it is enough to show that 
 \be\label{26}
S_{W_0}^\prime(0)(z) = c_1^0 \ne 0     
\end{equation} 
(in other words, that the zero set of the functional 
$J_1(f) = c_1$ is separated from the set of rotations (14) of the function $W_0$). This yields that the correspomdimg function (21), constructed by maximization of functional $J_1(f) = |c_1|$, is defined  and subharmonic on the whole disk $\D_{2|a_2^0|}$), and its maximaum is attained on the boundary circle.   

\bigskip
Assume, to the contrary, that $S_{W_0}^\prime(0)(z) = 0$. Then, by Lemma 2, 
$$
S_{W_0}(z) = \kappa_0 \circ \wh f_0(z) = c_0 + c_2 z^2 + 
c_3 z^3 + \dots \ , 
$$ 
where $\wh f_0$ is a holomorphic self-map of $\D$ of the form 
$$
\wh f_0(z) = \wh c_0 + \wh c_2 z^2 + \wh c_3 z^3 + \dots \ .  
$$
Since the function 
$$
\kappa_2(z) = \kappa_0(z^2) = 1/e + (2/e) z^2+ \dots
$$
also belongs to $\iota H_0^*$, it must be 
 \be\label{27} 
|c_2^0| > 2/e.   
\end{equation} 
Now consider the function 
$$ 
\wh f_1(z) = \sigma^{-1} \circ \Big\{\fc{\wh f_0 \circ  \sigma(z)}{z}\Big\} = \wh c_0 + \wh c_2 z + \wh c_3 z^2 + \dots \ ,  
$$ 
where
$$
\sigma(z) = (z - \wh c_0)/(1 - \ov{\wh c_0} z).
$$
This function also is a holomorphic self-map of the disk $\D$. Its  composition with $\kappa_0$ via (10), denoted by $f_1$, is a nonvanishing holomorphic self-map of $\D$, 
and a simple calculation, using (27), yields 
$$
f_1^\prime(0) = ( \kappa_0 \circ \wh f_1)^\prime(0) 
= |c_2^0| > 2/e, 
$$
which contradicts to Lemma 1. This proves the relations (25) and (26).

\bigskip\noindent
$\mathbf{8^0}$. \ Now we can finish the proof of the theorem. 

Take $n = 2$ and, letting $f_2(z) = f(z^2)$, consider on $H_0^\iy$ the plurisubharmonic functional 
  \be\label{28}
I_2(f) = \max \ (|J_2(f)|, |J_2(f_2)|). 
\end{equation} 
Similar to above, the lift of  this functional onto $\T_1$ generates via (19) a nonconstant radial
subharmonic function of on the disk (24). 
It is logarithmically convex, hence monotone increasing, and attains its maximal value at $|t| = 2 |a_2^0|$. 

By Parseval's equality for the boundary functions
$f(e^{i \theta}) = lim_{r\to 1} f(r e^{i \theta})$ of $f \in H_1^\infty$, we have   
$$
1 \ge \frac{1}{2 \pi} \int\limits_\pi^\pi |f(e^{i \theta})|^2 d \theta
= \sum_1^\infty |c_n|^2.  
$$ 
Applying it to the function 
$$
f(z) = \kappa_0(z) = \sum\limits_0^\iy c_n^0 z^n
$$
and noting that by (2), 
$|c_1^0|^2 = 4 e^{-2} = 0.541...$, one obtains that for this function, 
$$
\sum_2^\infty |c_n^0|^2 < 0.5 < |c_1^0|^2.  
$$
This implies (in view of the indicated connection of $|a_2^0|$ with $\kappa_0$) that the maximal value of the functional (28) on $H_0^\iy$ is attained on the functions $\kappa_0(z), \ \kappa_2(z) = \kappa_0(z^2)$ and equals 
$$
\max \ (|c_1^0|, |c_2^0|) = 2/e, 
$$
giving the desired estimate (1) for $n = 2$. 
The extremal extremal function is unique, up to rotations.

\bigskip
Now take $n = 3$ and, letting $f_3(z) = f(z^3)$, 
consider similar to (28) the functional 
$$
I_3(f) = \max \ (|J_3(f)|, |J_3(f_3)|)
$$
Arguing similar to the above case, one obtains 
$$ 
\max_{H_0^\iy} I_3(f) = \max \ (|c_1^0|, |c_3^0|) = 2/e, 
$$
giving the estimate (1) for $n = 3$. 

\bigskip
Taking subsequently $n = 4, 5, \dots$, one obtains by the same arguments that the estimate (1) is valid for all $n$, completing the proof of Theorem 1.

\newpage
\bigskip\bigskip
\centerline{\bf 4. PROOF OF THEOREM 2}

\bigskip
In view of the uniform convergence of $H^\iy$ functions 
on compact subsets of the unit disk, we have for the covering maps $\kappa_\rho: \ \D \to \mathcal A_\rho$, 
that their derivatives $\kappa_\rho^\prime(0)$ are  convergent to $\kappa_0^\prime(0) = 2/e < 1$ as $\rho \to 0$. 

Taking $\rho < \rho_0$ so small that $\kappa_\rho^\prime(0) < 1$, one can repeat the above 
arguments applied in the proof of Theorem 1 to the corresponding sets 
$H^0(\D, \mathcal A_\rho)$ with such $\rho$. 

But now 
we do not have the assertion on uniqueness of the covering map $\kappa_\rho$ on which the maximal value (3) is attained.

\bigskip\bigskip
\centerline{\bf 5. REMARK ON THE HUMMEL-SCHEINBERG-ZALCMAN
CONJECTURE}

\bigskip
The Krzyz conjecture was extended in 1977 by Hummel, Scheinberg and
Zalcman to arbitrary Hardy spaces $H^p, \ p > 1$ on the unit disk,
for which there is conjectured that the coefficients of nonvanishing
functions $f(z) \in H^p, \ p > 1$, with $\|f\|_p \le 1$ satisfy
$$
|c_n| \le (2/e)^{1-1/p},
$$
with equality for the function
$$
f_n(z) = \Bigl[\frac{(1 + z^n)^2}{2}\Bigr]^{1/p} \ \Bigl[\exp
\frac{z^n - 1}{z^n + 1}\Bigr]^{1-1/p}
$$
and its rotations (see \cite{HSZ}). As $p \to \infty$, this yields Krzyz's conjecture for $H^\infty$ (without
uniqueness of extremal functions)

This problem also has been investigated by many authors, but it
still remains open. The only known results here are that the
conjecture is true for $n = 1$ proved by Brown \cite{Br} as well as
some results for special subclasses of $H^p$, see \cite{Br},
\cite{BW}, \cite{Su}.

Some important intrinsic features of $H^\infty$ functions,
essentially involved in the proof of Krzyz's conjecture, are lost in
$H^p$. However, the above arguments can be appropriately modified and completed to include also the Hummel-Scheinberg-Zalcman conjecture. This will be presented in a separate work.

\bigskip
{\small\em{ \leftline{Department of Mathematics, Bar-Ilan
University, Ramat-Gan, Israel} \leftline {Department of Mathematics,
University of Virginia, Charlottesville, VA 22904-4137, USA}}}


\bigskip
\bigskip
\begin{thebibliography}{MSTV}
{\small


\bibitem{AMC}
J. Agler and J. McCarthy, {\it The Krzyz conjecture and an Entropy
conjecture}, arXiv:1803.09718v1 [mathCV].

\bibitem{Ba}
R. W. Barnard, {\it Open problems and conjectures in complex
analysis}, Computational Methods and Function Theory Proceedings,
Valparaiso 1989, St. Ruscheweyh, E. B. Saff, L. C. Salinas, R. S.
Varga (eds.), Lecture Notes in Mathematics 1435, Springer, Berlin,
1990, pp. 1-26.

\bibitem{Be}
L. Bers, {\it Fiber spaces over Teichm\"{u}ller spaces}, Acta Math.
\textbf{130} (1973), 89-126.

\bibitem{Br}
J. E. Brown, {\it A coefficient problem for nonvanishing $H^p$
functions}, Complex Variables \textbf{4} (1985), 253-265.

\bibitem{BW}
J. E. Brown and J. Walker, {\it A coefficient estimate for
nonvanishing $H^p$ functions}, Rocky Mountain J. Math. \textbf{18}
(1988), 707-718.

\bibitem{Ea}
C.J. Earle, {\it On quasiconformal extensions of the
Beurling-Ahlfors type}, Contribution to Analysis, Academic Press,
New York, 1974, pp. 99-105.

\bibitem{EKK}
C.J. Earle, I. Kra  and S.L. Krushkal, {\em Holomorphic motions and Teichm\"{u}ller spaces}, Trans. Amer. Math. Soc. \textbf{944} (1994), 927-948. 

\bibitem{GL}
F.P. Gardiner and N. Lakic, {\it Quasiconformal \Te \ Theory}, Amer.
Math. Soc., 2000.

\bibitem{Go}
G.M. Goluzin, {\it Geometric Theory of Functions of Complex
Variables}, Transl. of Math. Monographs, vol. 26, Amer. Math. Soc.,
Providence, RI, 1969.

\bibitem{Ho}
Ch. Horowitz, {\it Coefficients of nonvanishing functions in
$H^\iy$}, Israel J. Math. \textbf{30} (1978), 285-291.

\bibitem{HSZ}
J. A. Hummel, S. Scheinberg and L. Zalcman, {\it A coefficient
problem for bounded nonvanishing functions}, J. Anal. Math.
\textbf{31} (1977), 169-190.

\bibitem{Kr1}
S.L. Krushkal, {\it Hyperbolic geodesics, Krzyz's conjecture and
beyond}, arXiv:1603.02668[mathCV].

\bibitem{Kr2}
S.L. Krushkal, {\it A general coefficient theorem for univalent
functions}, arXIv:1908.05183[mathCV].

\bibitem{Kz}
J. Krzyz, {\it Coefficient problem for bounded nonvanishing
functions}, Ann. Polon. Math. \textbf{70} (1968), 314.

\bibitem{LS}
Z. Lewandowski and J. Szynal, {\it On the Krzyz conjecture and
related problems}, XVIth Rolf Nevanlinna Colloquium (Joensu, 1995),
de Gruyter, Berlin, 1996, pp. 257-269.

\bibitem{MSUV}
M.J. Martin,, E.T.Sawyer, I. Uriatre-Tuero and D. Vukovich, {\it The
Krzyz conjecture revisited}, arXiv:1311.7668 [mathCV]

\bibitem{Ma}
W.S. Massey, {\it A Basic Course in Algebric Topology}, Springer,
New York, 1991.

\bibitem{PS}
D. V. Prokhorov and J. Szynal, {\it Coefficient estimates for
bounded nonvanishing functions}, Bull. Acad. Polon. Sci. Ser.
\textbf{29} (1981), 223-230.

\bibitem{R}
W.W. Rogosinski, {\it On the order of the derivatives of a function
analytic in an angle}, J. London Math. Soc. \textbf{20} (1945),
100-109.

\bibitem{Sa}
N. Samaris, {\it A proof of Krzyz's conjecture for the fifth
coefficient}, Complex Variables \textbf{48} (2003), 753-766.

\bibitem{Su}
T.J. Suffridge, {\it Extremal problems for nonvanishing $H^p$
functions}, Computational methods and function theory (Valparaiso,
1989), Lectures Notes in Math., 1435, Springer, Berlin, 1990, pp.
177-190.

\bibitem{Sz}
W. Szapiel, {\it A new approach to the Krzyz conjecture}, Ann. Univ.
Mariae Curie-Sklodowska Sect. A \textbf{48} (1994), 169-192.

\bibitem{Ta}
D.-L. Tan, {\it Estimates of coefficients of bounded nonvanishing
analytic functions}, Chinese Ann. Math. Ser. A \textbf{l4} (1983)
97-104.

\bibitem{Te}
O. Teichm\"{u}ller. {\it Ein Verschiebungssatz der quasikonformen
Abbildung}, Deutsche Math. \textbf{7} (1944), 336-343.

\bibitem{Ve}
E. Vesentini, {\em Complex geodesc and holomorphic mappings},
Sympos. Math. \textbf{26} (1982), 211-230.



}

\end{thebibliography}
\end{document}